\magnification=\magstep 1
\input amstex
\documentstyle{amsppt}

\def\O{\Omega}
\def\D{\Delta}

\def\n{\nabla}

\def\RR{\Bbb R}

\def\b{\beta}
\def\Si{\Sigma}

\def\CD{{\Cal D}_{\rho}}
\def\a{\alpha}
\def\a1{\alpha_{\rho}}
\def\e{\epsilon}
\def\d{\delta}

\topmatter
\title Existence results for mean field equations
\endtitle
\author Weiyue Ding, J\"urgen Jost, Jiayu Li and Guofang Wang
\endauthor
\address Institute of Mathematics, Academia Sinica, Beijing 
100080, P. R. China
\endaddress 
\email dingwy\@public.intercom.co.cn
\endemail
\address Max-Planck-Institute for Mathematics in the Sciences, 
Inselstr. 22-26, 04103 Leipzig, Germany
\endaddress
\email jost\@mis.mpg.de
\endemail
\address Institute of Mathematics, Academia Sinica, 100080 Beijing,
China
\endaddress
\email lijia\@public.intercom.co.cn
\endemail
\address Institute of Mathematics, Academia Sinica, 100080 Beijing,
China
\endaddress
\email gwang\@mis.mpg.de
\endemail
%\keywords 
%\endkeywords
%\subjclass 
%\endsubjclass
\abstract 
Let $\O$ be an annulus. We prove that the mean field equation
$$\eqalign{
-\D\psi&=\frac{e^{-\b\psi}}{\int_{\O}e^{-\b\psi}} \quad  \hbox{in }\O\cr
\psi &=0  \qquad\qquad \hbox{ on }\partial\O\cr}
$$
admits a solution for $\b\in (-16\pi,-8\pi)$. This is a supercritical case for
the Moser-Trudinger inequality.
 \endabstract
\endtopmatter

\document
\head 1.Introduction
\endhead

Let $\O$ be a smooth bounded domain in $\RR^2$. In this paper, we consider 
the following mean field equation
$$
\eqalign{
-\D \psi&= \frac{e^{-\b\psi}}{\int_{\O}e^{-\b\psi}}, \quad\hbox{in } \O,\cr
\psi&=0,\qquad\qquad\hbox{ on }\partial\O,\cr}\leqno(1.1)
$$
for $\b\in (-\infty, +\infty)$. (1.1) is the Euler-Lagrange equation of the
following functional
$$J_{\b}(\psi)=\frac12\int_{\O}|\n \psi|^2+\frac1{\b}\log\int_{\O}e^{-\b\psi}
\leqno(1.2)$$
in $H^{1,2}_0(\O)$. This variational problem arises from Onsager's 
vortex model for turbulent Euler flows. In that interpretation, $\psi$ is the
stream function in the infinite vortex limit, see [MP,p256ff]. The corresponding
canonical Gibbs measure and partition function are finite precisely
if $\b>-8\pi$. In that situation, Caglioti et al. [CLMP1]
and Kiessling [K] showed the existence of a minimizer of $J_{\b}$. 
This is based on the Moser-Trudinger inequality
$$\frac12\int_{\O}|\n\psi|^2\ge \frac1{8\pi}\log\int_{\O}e^{-8\pi\psi},
\qquad\hbox{for any }\psi\in H^{1,2}_0(\O),\leqno(1.3)$$
which implies the relevant compactness and coercivity condition for $J_{\b}$
in case $\b>-8\pi$. For $\b\le -8\pi$, the situation becomes different 
as described in [CLMP1]. On the unit disk, solutions blow up if one approaches 
$\b=-8\pi$ -the critical case for (1.3)-(see also [CLMP2] and [Su]),
and more generally, 
on starshaped domains, the Pohozaev identity yields a lower bound on the 
possible values of $\b$ for which solutions exist. On the other hand, for an
annulus, [CLMP1] constructed radially symmetric solutions for any $\b$, and
the construction of Bahri-Coron [BC] makes it plausible that solutions on 
domains with non-trivial topology exist below $-8\pi$. Thus, for $\beta\le
-8\pi$, $J_{\b}$ is no longer compact and coercive in general, and the existence
of solution depends on the geometry of the domain.

In the present paper, we thus consider the supercritical case $\b<-8\pi$ on 
domains with non-trivial topology.

\proclaim{Theorem 1.1}Let $\O\subset \RR^2$ be a smooth, bounded domain
whose complement contains a bounded region,  e.g. $\O$ an annulus. Then
(1.1) has a solution for all $\b\in (-16\pi,-8\pi)$.
\endproclaim

The solutions we find, however, are not minimizers of $J_{\b}$-those 
do not exist in case $\b<8\pi$, since $J_{\b}$ has no lower bound-but 
unstable
critical points. Thus, these solutions might not be relevant to the 
turbulence
problem that was at the basis of [CLMP1] and [K].

Certainly we can generalize Theorem 1.1 to the following equation
$$
\eqalign{
-\D \psi&= \frac{Ke^{-\b \psi}}{\int_{\O}Ke^{-\b \psi}}, \quad\hbox{in } \O,\cr
\psi&=0,\qquad\qquad\hbox{ on }\partial\O,\cr}
$$
which was  studied in [CLMP2]. Here $K$ is a positive function on $\bar{\O}$.

With the same method, we may also handle the equation
$$\D u-c+cKe^u=0,\qquad \hbox{for } 0\le c<\infty\leqno(1.4)$$
on a compact Riemann surface $\Si$ of genus at least $1$, where $K$ is a
positive function

(1.4) can also be considered as a mean field equation because it is the
Euler-Lagrange equation of the functional
$$J_c(u)=\frac12\int_{\Si}|\n u|^2+c\int_{\Si}u-c\log\int_{\Si}Ke^u.
\leqno(1.5)$$
Because of the term $c\int_{\Si}u$, $J_c$ remains invariant under adding 
a constant to $u$, and therefore we may normalize $u$ by the condition
$$\int_{\Si}Ke^u=1$$
which explains the absence of the factor $(\int Ke^u)^{-1}$
in (1.4). $c<8\pi$ again is a subcritical case that can easily be handled 
with the 
Moser-Trudinger inequality. The critical case $c=8\pi$ yields the so-called
Kazdan-Warner equation [KW] and was treated in [DJLW] and [NT] 
by giving sufficient
conditions for the existence of a minimizer of $J_{8\pi}$. Here, we construct 
again saddle point type critical points to show
\proclaim{Theorem 1.2} Let $\Si$ be a compact Riemann surface of positive
genus. Then (1.4) admits a non-minimal solution for $8\pi< c<16\pi.$
\endproclaim

Now we give a outline of the proof of the Theorems. First from the
non-trivial topology of the domain, we can define a minimax value 
$\alpha_{\b}$, which is bounded below by an improved Moser-Trudinger 
inequality, for $\b\in(-16\pi,-8\pi)$. Using a trick introduced by Struwe
in [St1] and [St2], for a certain dense subset 
$\Lambda\subset (-16\pi,-8\pi)$ we can overcome the lack of a
coercivity condition and show that $\alpha_{\b}$ is achieved by some 
$u_{\b}$ for $\b\in\Lambda$. Next, for any fixed 
$\bar{\b}\in (-16\pi,-8\pi)$, considering a sequence $\b_k\subset\Lambda$ 
tending to 
$\bar{\b}$, with the help of results in [BM] and [LS] 
we show that $u_{\b_k}$ subconverges strongly to some $u_{\bar{\b}}$
which achieves $\alpha_{\bar{\b}}$.

After completing our paper, we were informed that Struwe and Tarantello
[ST] obtained a non-constant solution of (1.4), when $\Si$ is a flat torus with 
fundamental cell domain $[-\frac12, \frac12]\times [-\frac12, \frac12]$,
$K\equiv 1$ and $c\in (8\pi, 4\pi^2)$. In this case, it is easy to
check that our solution obtained in Theorem 1.2 is non-constant.

Our research was carried out at the Max-Planck-Institute for 
Mathematics in the Sciences in Leipzig. The first author thanks the 
Max-Planck-Institute for the hospitality and good working conditions.
The third author was supported
by a fellowship of the Humboldt foundation, whereas the fourth author
was supported by the DFG through the Leibniz award of the second author.

\head 2.  Minimax values
\endhead
 Let $\rho=-\beta$ and $u=-\beta\psi$. We rewrite (1.1) as
$$
\eqalign{
-\Delta u&=\rho\frac{e^u}{\int_{\O}e^u}, \quad \hbox{ in }\O,\cr
u&=0,\qquad\qquad \hbox{on } \partial\O,\cr}\leqno (2.1)
$$
and (1.2) as
$$J_{\rho}(u)=\frac12\int_{\O}|\n u|^2-\rho\log\int_{\O}e^u\leqno (2.2)$$
for $u\in H^{1,2}_0(\O).$

It is easy to see that $J_{\rho}$ has no lower bound for $\rho\in(8\pi,16\pi)$.
Hence, to get a solution of (1.1) for $\rho\in(8\pi,16\pi)$, we have to use
a minimax method. First, we define a center of mass of $u$ by
$${m_c}(u)=\frac{\int_{\O}xe^u}{\int_{\O}e^u}.$$
Let $B$ be the bounded component of $\RR^2\setminus \O$. For simplicity,
we assume that $B$ is the unit disk centered at the origin.
Then we define a family of functions
$$h:D\to H^{1,2}_0(\O)$$
satisfying
$$\lim_{r\to 1}J_{\rho}(h(r,\theta))\to -\infty\leqno(2.3)$$
and
$$ \lim_{r\to 1}{m_c}(h(r,\theta)) \hbox{ is a continuous curve enclosing 
$B$.}\leqno(2.4)$$
Here $D=\{(r,\theta)|0\le r< 1,\theta\in [0,2\pi)\}$ is the open unit disk.
We denote the set of all such families by $\CD$.
It is easy to check that $\CD\not = \emptyset$.
Now we can define a minimax value
$$\a1 :=\inf_{h\in\CD}\sup_{u\in h(D)}J_{\rho}(u).$$
The following lemma will make crucial use of the non-trivial topology of
$\O$, more precisely of the fact that the complement of $\O$ has a bounded
component.
\proclaim{Lemma 2.1} $\a1>-\infty$ for any $\rho \in (8\pi,16\pi)$.
\endproclaim
\demo{Remark} It is an interesting question weather $\alpha _{16\pi}=-\infty.$
\enddemo

To prove Lemma 2.1, we use the improved Moser-Trudinger inequality of [CL]
(see also [A]). Here we have to modify a little bit.
\proclaim{Lemma 2.2} Let $S_1$ and $S_2$ be two subsets of $\bar{\O}$ satisfying
$dist(S_1,S_2)\ge\delta_0>0$ and $\gamma_0\in (0,1/2)$.
For any $\e>0$, there exists a constant $c=c(\e,\delta_0,\gamma_0)>0$
such that
$$\int_{\O}e^u\le c\exp\{\frac1{32\pi-\e}\int_{\O}|\n u|^2+c\}$$
holds for all $u\in H^{1,2}_0(\O)$ satisfying
$$\frac{\int_{S_1}e^u}{\int_{\O}e^u}\ge \gamma_0\qquad \hbox{and}\qquad
\frac{\int_{S_2}e^u}{\int_{\O}e^u}\ge \gamma_0.\leqno(2.5)$$
\endproclaim
\demo{Proof}The Lemma follows from the argument in [CL] and the following
Moser-Trudinger inequality
$$\frac12\int_{\O}|\n u|^2-8\pi\log\int_{\O}e^u\ge c\leqno(*)$$
for any $u\in H_0^{1,2}(\O)$, where $c$ is a constant independent of 
$u\in H_0^{1,2}(\O)$.\qed
\enddemo
We will discuss the inequality $(*)$ and its application in another paper.
\demo{Proof of Lemma 2.1} For fixed $\rho\in (8\pi,16\pi)$ we claim that 
there exists a constant $c_{\rho}$ such that
$$\sup_{u\in h(D)}J_{\rho}(u)\ge c_{\rho}, \qquad \hbox{for any } h\in \CD.
\leqno (2.6)$$
Clearly (2.6) implies the Lemma. By the definition of $h$, for any $h\in \CD$,
there exists $u\in h(D)$ such that 
$${m_c}(u)=0.$$
We choose $\e>0$ so small that $\rho< 16\pi-2\e$. Assume (2.6) does not hold. 
Then we have sequences $\{h_i\}\subset\CD$ and $\{u_i\}\subset H^{1,2}_0(\O)$
such that $u_i\in h_i(D)$ and
$${m_c}(u_i)=0\leqno(2.7)$$
$$\lim_{i\to\infty}J(u_i)=-\infty.\leqno(2.8)$$

%For simplicity, we assume that $\\int_{\O}e^{u_i}=1$. (Recall that
%$J_{\rho}$ is invariant 
We have the following Lemma.
\proclaim{Lemma 2.3} There exists $x_0\in \bar{\O}$ such that
$$\lim_{i\to \infty}
\frac{\int_{B_{1/2}(x_0)\cap \O}e^{u_i}}{\int_{ \O}e^{u_i}}\to 1.
\leqno(2.9)$$
\endproclaim
\demo{Proof} Set 
$$A(x):=
\lim_{i\to\infty}\frac{\int_{B_{1/4}(x)\cap \O}e^{u_i}}{\int_{ \O}e^{u_i}}.
$$
Assume that the Lemma were false,
then there exists $x_0\in \bar{\O}$ such that
$$A(x_0)<1 \quad \hbox{and}\quad A(x_0)\ge A(x) \quad\hbox{for any } x\in\O.$$
It is easy to check $A(x_0)>0$, since $\O$ can be covered by finite many
balls of radius $1/4$. Let $\gamma_0=A(x_0)/2$.
 Recalling (2.8) and applying lemma 2.2, we obtain
$$\frac{\int_{\O\setminus B_{1/2}(x_0)}e^{u_i}}
{\int_{\O}e^{u_i}}\to 0\leqno(2.10)$$
as $i\to \infty$, which implies (2.9). \qed
\enddemo
Now we continue to prove Lemma 2.1. (2.9) implies
$$
\eqalign{
\frac{\int_{\O}xe^{u_i}}{\int_{\O}e^{u_i}}-x_0&
=\frac{\int_{\O}(x-x_0)e^{u_i}}{\int_{\O}e^{u_i}}\cr
&=\frac{\int_{B_{1/2}(x_0}(x-x_0)e^{u_i}}{\int_{\O}e^{u_i}}+o(1)\cr}
$$
which, in turn, implies  that $|{m_c}(u_i)-x_0|<2/3$. 
This contradicts (2.7).\qed
\enddemo

%\demo{Remark 2.4}
%For $\rho\in (0,8\pi]$, we can also define a minimax value of $J_{\rho}$
%by changing (2.3) as follows
%$$\lim_{r\to 1}J_{\rho}(h(r,\theta))\to 
%\inf_{u\in H^{1,2}_0(\O)} J_{\rho}.\leqno(2.3)'$$
%(From the Moser-Trudinger inequality (1.3),  
%$\inf_{u\in H^{1,2}_0(\O)}J_{\rho}$ 
%is finite.) By the same method given above, we can prove that
%$$\a1>\inf_{u\in H^{1,2}_0(\O)} J_{\rho}\leqno(2.11)$$
%for $\rho\in (0,8\pi]$. For $\rho\in (0,8\pi)$, as mentioned in the
%introduction, $J_{\rho}$ satisfies a compactness and coercivity condition.
%Hence, it is easy to show that $\a1$ is achieved by some $u_{\rho}$
%which is not a minimizer of $J_{\rho}$. For $\rho=8\pi$,
%see Remark 4.2 below.
%\enddemo

\proclaim{Lemma 2.4} $\a1/\rho$ is non-increasing in $(8\pi,16\pi)$.
\endproclaim
\demo{Proof} We first observe that if $J(u)\le 0$, then $\log\int_{\O}e^u>0$ 
which implies that
$$J_{\rho}(u)\ge J_{\rho'}(u) \qquad \hbox{for }\rho'\ge\rho.$$
Hence ${\Cal D}_{\rho}\subset {\Cal D}_{\rho'}$ for any $16\pi>
\rho'\ge\rho>8\pi$. On the other hand, it is clear that
$$\frac{J_{\rho}}{\rho}-\frac{J_{\rho'}}{\rho'}=\frac12(\frac1{\rho}-
\frac1{\rho'})\int_{\O}|\n u|^2\ge 0,$$
if $\rho'\ge \rho$. Hence we have
$$\frac{\a1}{\rho}\ge \frac{\alpha_{\rho'}}{\rho'}$$
for $16\pi>\rho'\ge\rho>8\pi$. \qed
\enddemo

\head 3. Existence for a dense set
\endhead

In this section we show that $\a1$ is achieved if $\rho$ belongs to a 
certain
dense subset of $(8\pi,16\pi)$ defined below.

The crucial problem for our functional is the lack of a coercivity 
condition, 
{\it i.e.} for a Palais-Smale sequence $u_i$ for $J_{\rho}$, we do not know
whether $\int_{\O}|\n u_i|^2$ is bounded.

We first have the following lemma.
\proclaim{Lemma 3.1} Let $u_i$ be a Palais-Smale sequence for $J_{\rho}$, 
i.e.
$u_i$ satisfies
 $$|J_{\rho}(u_i)|\le c<\infty\leqno(3.1)$$
 and
 $$dJ_{\rho}(u_i)\to 0 \hbox{ strongly in } H^{-1,2}(\O)\leqno(3.2).$$

If, in addition, we have
$$\int_{\O}|\n u_i|^2\le c_0, \qquad\hbox{for } i=1,2,\cdots\leqno(3.3)$$
for a constant $c_0$ independent of $i$, then $u_i$ subconverges 
to a critical point $u_0$ for $J_{\rho}$ strongly in $H^{1,2}_0(\O)$.
\endproclaim
\demo{Proof} The proof is standard, but we provide it here for convenience 
of the reader.

Since $\int_{\O}|\n u_i|^2$ is bounded, there exists $u_0\in H^{1,2}_0(\O)$ 
such that
\roster
\item"(i)"$u_i$ converges to $u_0$ weakly in $H^{1,2}_0(\O),$
\item"(ii)"$u_i$ converges to $u_0$ strongly in $L^p(\O)$ for any $p>1$ and
almost everywhere,
\item"(iii)"$e^{u_i}$ converges to $e^{u_0}$ strongly in $L^p(\O)$ for 
any $p\ge 1$.
\endroster
From (i)-(iii), we can show that $dJ(u_0)=0$, {\it i.e.} $u_0$ satisfies
$$-\D u_0=\rho \frac{e^{u_0}}{\int_{\O}{e^{u_0}}}.$$
Testing $dJ_{\rho}$ with $u_i-u_0$, we obtain
$$\eqalign{
o(1)&=\langle dJ_{\rho}(u_i)-dJ_{\rho}(u),u_i-u_0\rangle\cr
&=\int_{\O}|\n(u_i-u_0)|^2-\rho\int_{\O}(\frac{e^{u_i}}{\int_{\O}e^{u_i}}-
\frac{e^{u_0}}{\int_{\O}e^{u_0}})(u_i-u_0)\cr
&=\int_{\O}|\n(u_i-u_0)|^2+o(1),\cr
}$$
by (i)-(iii).
Hence $u_i$ converges to $u_0$ strongly in $H^{1,2}_0(\O)$.
\qed
\enddemo

Since by Lemma 2.4 $\rho\to\a1/\rho$ is non-increasing in $(8\pi,16\pi)$,
$\rho\to\a1/\rho$ is a.e. differentiable. Set
$$\Lambda:=\{\rho\in(8\pi,16\pi)|\a1/\rho \hbox{ is differentiable at }\rho\}.
\leqno(3.4)$$
$\bar\Lambda=[8\pi,16\pi]$, see  [St1]. Let $\rho\in\Lambda$ and choose 
$\rho_k\nearrow\rho$ such that 
$$0\le\lim_{k\to\infty}-\frac1{(\rho-\rho_k)}(\frac{\a1}{\rho}-\frac{\alpha_{\rho_k}}{\rho_k})
\le c_1\leqno(3.5)$$
for some constant $c_1$ independent of $k$.
\proclaim{Lemma 3.2} $\a1$ is achieved by a critical point $u_{\rho}$ for
$J_{\rho}$ provided that $\rho\in \Lambda$.
\endproclaim
\demo{Proof} Assume, by contradiction, that
the Lemma were false. From Lemma 3.1, there exists $\delta>0$ such that
$$\|dJ_{\rho}(u)\|_{H^{-1,2}(\O)}\ge 2\delta\leqno(3.6)$$
in
$$
N_{\delta}:=\{u\in H^{1,2}_0(\O)|\int_{\O}|\n u|^2\le c_2, |J_{\rho}(u)-\a1|
<\d\}.$$
Here, $c_2$  is any fixed constant such that $N_{\d}\not =\emptyset$.
Let $X_{\rho}:N_{\d}\to H^{1,2}_0(\O)$ be a pseudo-gradient vector
field for $J_{\rho}$ in $N_{\d}$, {\it i.e.} a locally Lipschitz vector
field of norm $\|X_{\rho}\|_{H^{1,2}_0}\le 1$ with
$$\langle dJ_{\rho}(u), X_{\rho}(u)\rangle < -\d.\leqno(3.7)$$
See [P]  for the construction of $X_{\rho}$.

Since 
$$\eqalign{
\|dJ_{\rho}(u)-dJ_{\rho_k}(u)\|&= \|dJ_{\rho}-\frac{\rho}
{\rho_k}dJ_{\rho_k}(u)\|
+\|(1-\frac{\rho}{\rho_k})dJ_{\rho_k}(u)\|\cr
&\le \frac12(1-\frac{\rho}{\rho_k})\int|\n u|^2+
c(1-\frac{\rho}{\rho_k})\int_{\O}|\n u|^2\to 0\cr}
$$
uniformly in $\{u|\int_{\O}|\n u|^2\le c_2\}$, $X_{\rho}$ is also a 
pseudo-gradient vector field for $J_{\rho_k}$ in $N_{\d}$ with
$$\langle dJ_{\rho_k}(u), X_{\rho}(u)\rangle < -\d/2,\leqno(3.8)$$
for $u\in N_{\d}$, provided that $k$ is sufficiently large.

For any sequence $\{h_k\}$,  $h_k\in {\Cal D}_{\rho_k}\subset \CD$ such that
$$\sup_{u\in h_k(D)}J_{\rho_k}(u)\le \alpha_{\rho_k}+\rho-\rho_k\leqno(3.9)$$
and all $u\in h_k(D)$ such that 
$$J_{\rho}(u)\ge \a1-(\rho-\rho_k),\leqno(3.10)
$$
we have the following estimate
$$\eqalign{
\frac12\int_{\O}|\n u|^2&=\rho\cdot\rho_k 
\frac {\frac{J_{\rho_k}(u)}{\rho_k}-\frac{J_{\rho}(u)}{\rho}}{\rho-\rho_k}
\cr
&\le \rho\cdot\rho_k\frac 
{\frac{\alpha_{\rho_k}}{\rho_k}-\frac{\alpha_{\rho}}{\rho}}{\rho-\rho_k}
+(\rho+\rho_k)\cr
&\le C\cr}\leqno(3.11)$$
by (3.5), (3.9) and (3.10), where $C=(16\pi)^2c_1+32\pi.$

Now we consider in $N_{\d}$ the following pseudo-gradient flow for $J_{\rho}$.
First choose a Lipschitz continuous cut-off function $\eta$ such that 
$0\le\eta\le 1$, $\eta=0$ outside $N_{\d}$, $\eta=1$ in $N_{\d/2}$. Then
consider the following flow in $H^{1,2}_0(\O)$ generated by $\eta X_{\rho}$
$$\eqalign{
\frac{\partial \phi}
{\partial t}(u,t)&=\eta(\phi(u,t))X_{\rho}(\phi(u,t))\cr
\phi(u,0)&=u.\cr}
$$
By (3.7) and (3.8), for $u\in N_{\d/2}$, we have
$$\frac d{dt}J_{\rho}(\phi(u,t))_{|_{t=0}}\le-\d \leqno(3.12)
$$
and
$$\frac d{dt}J_{\rho_k}(\phi(u,t))_{|_{t=0}}\le-\d/2 \leqno(3.13)
$$
for large $k$. 

It is clear that for any $h\in {\Cal D}_{\rho_k}$ $h(r,\theta)\not\in N_{\d}$
for $r$ close to $1$. Hence $\phi(h,t)\in {\Cal D}_{\rho_k}$ for any $t>0$.
In particular, $\phi(\cdot,t)$ preserves the class of $h_k\in 
{\Cal D}_{\rho_k}$ with condition (3.9). On the other hand, for any $h\in 
{\Cal D}_{\rho}$ by definition
$$\sup_{u\in h(D)}J_{\rho}(u)\ge \a1.$$
Hence for any $h_k\in {\Cal D}_{\rho_k}$ with condition (3.9), 
$\sup_{u\in\phi(h(D),t)}J_{\rho}(u)$ is achieved in $N_{\d/2}$, provided that $k$
is large enough. Consequently, by (3.12), we have
$$\frac d{dt}\sup\{J_{\rho}(u)|u\in \phi(h(D),t)\}\le -\d
$$ for all $t\ge 0$, which is a contradiction.
\qed
\enddemo

\head 4. Proof of Theorem 1.1
\endhead
From section 3, we know that for any $\bar{\rho}\in (8\pi,16\pi)$ there exists
a sequence $\rho_k\nearrow\bar{\rho}$ such that $\alpha_{\rho_k}$ is achieved by 
$u_k$. Consequently $u_k$ satisfies
$$
\eqalign{
-\D u_k&=\rho_k\frac {e^{u_k}}{\int_{\O}e^{u_i}}, \quad\hbox{in } \O,\cr
u_k&=0,\qquad\qquad\hbox{ on }\partial \O.\cr}\leqno(4.1)$$
From Lemma 2.4, we have 
$$J_{\bar{\rho}}(u_k) =\alpha_{\rho_k}\hbox{ is bounded.}\leqno(4.2)$$
%It, in turn, implies that 
%$$\int_{\O}e^{u_k}\ge c_0\leqno(4.2)$$
for some constant $c_0>0$ which is independent of $k$.
Let $v_k=u_k-\log\int_{\O}e^{u_k}$. Then $v_k$ satisfies
$$-\D v_k=\rho_k e^{v_k}\leqno(4.3)$$
with
$$\int_{\O}e^{v_k}=1.\leqno(4.4)$$

By  results of Brezis-Merle [BM] and Li-Shafir [LS] we have
\proclaim{Lemma 4.1} ([BM], [LS]) There exists a subsequence (also denoted
by $v_k$) satisfying one of the following alternatives:
\roster
\item"(i)"$\{v_k\}$ is bounded in $L^{\infty}_{loc}(\O)$;
\item"(ii)"$v_k\to -\infty$ uniformly on  any compact subset of $\O$; 
\item"(iii)" there exists a finite blow-up set $\Si=\{a_1,\cdots,a_m\}
\subset
\O$ such that, for any $1\le i\le m$, there exists $\{x_k\}\subset\O$, $
x_k\to a_i$, $u_k(x_k)\to\infty$, and $v_k(x)\to -\infty$ uniformly on any
compact subset of $\O\setminus \Si$. Moreover,
$$\rho_k\int_{\O}e^{v_k}\to \sum_{i=1}^m8\pi n_i\leqno(4.5)$$
where $n_i$ is positive integer.
\endroster
\endproclaim

For our special functions $v_k$, 
we can improve Lemma 4.1 as follows

\proclaim{Lemma 4.2}  There exists a subsequence (also denoted
by $v_k$) satisfying one of the following alternatives:
\roster
\item"(i)"$\{v_k\}$ is bounded in $L_{loc}^{\infty}(\O)$;
\item"(ii)"$v_k\to -\infty$ uniformly on   $\bar{\O}$; 
\item"(iii)" there exists a finite blow-up set $\Si=\{a_1,\cdots,a_m\}
\subset
\bar{\O}$ such that, for any $1\le i\le m$, there exists $\{x_k\}\subset\O$, $
x_k\to a_i$, $u_k(x_k)\to\infty$, and $v_k(x)\to -\infty$ uniformly on any
compact subset of $\bar{\O}\setminus \Si$. Moreover, (4.5) holds.
\endroster
\endproclaim
\demo{Proof}From Lemma 4.1, we only have to consider one more case in which
 blow-up points  are in the boundary of $\O$. There are two possibilities:
 One is  bubbling too fast such that after rescaling we obtain a
 solution of $\D u=e^u$ in a half plane; Another is bubbling slow
 such that after rescaling we obtain a solution 
of $\D u=e^u$ in $\RR^2$. One can exclude the first case. In the second case,
one can follow the idea in [LS] to show that (4.5) holds. See also [L].\qed
\enddemo 

\demo{Proof of Theorem 1.1} 
 (4.4), (4.5) and $\bar{\rho}\in (8\pi,16\pi)$ 
imply that cases (ii) and (iii) in Lemma 4.2   does not occur. 
%It is clear that
%$u_k \ge 0$ by the maximum principle. Hence from (4.3), we have
%$v_k\ge -\log c_0$ which excludes (ii). 
Consequently $\{v_k\}$ is bounded in $L^{\infty}_{loc}(\O)$. 
Now we can again apply Lemma 2.2 as follows.

Let $S_1$ and $S_2$ be two disjoint compact subdomains of $\O$. Since 
$\{v_k\}$ is bounded in $L^{\infty}_{loc}(\O)$,
we have
$$\frac{\int_{S_i}e^{u_k}}{\int_{\O}e^{u_k}}=\int_{S_i}e^{v_k}\ge c_0, \qquad
i=1,2$$
for a constant $c_0=c_0(S_1,S_2,\O)>0$ independent of $k$. Choosing $\e$
such that $16\pi-\bar{\rho}>2\e$ and applying Lemma 2.2, with the help of (4.2),
we obtain 
$$\eqalign{
c&\ge J_{\rho_k}(u_k)=\frac12\int_{\O}|\n u_k|^2-{\rho_k}
\log\int_{\O}e^{u_k}\cr
&\ge\frac12(1-\frac{{\rho_k}}{16\pi-\e/2})\int_{\O}|\n u|^2\cr
&\ge\frac12(1-\frac{\bar{\rho}}{16\pi-\e/2})\int_{\O}|\n u|^2\cr
}$$
which implies
that $\int_{\O}|\n u_k|^2$ is bounded. Now by the same argument in the proof
of Lemma 3.1,
 $u_k$ subconverges to $u_{\bar{\rho}}$ strongly in $H^{1,2}_0(\O)$ and
$u_{\bar{\rho}}$ is a critical point of $J_{\bar{\rho}}$. Clearly,
$u_{\bar{\rho}}$ achieves $\alpha_{\bar{\rho}}$. This finishes the
proof of  Theorem 1.1.

\qed
\enddemo
%\demo{Remark 4.2} Recall that for $\rho=8\pi$, we can also define 
%a minimax value satisfying (2.11), see Remark 2.4. By the argument presented
%above, we only need  to exclude case (iii) in Lemma 4.1 with $m=1$ and $n_1=1$.
%In this  case, we can show that $\alpha_{8\pi}=\inf_{u\in H^{1,2}_0(\O)}
%J_{8\pi}(u)$ which contradicts (2.11).
%\enddemo

%Therefore, Theorem 1.1, Remarks 2.4 and 4.2 imply that (1.1) has a solution,
%which is not a minimizer of $J_{\beta}$,
%for any $\beta\in (-16\pi,0)$.

\demo{Proof of Theorem 1.2} Since the proof is very similar to one presented above,
 we only give a
 sketch of the proof of Theorem 1.2. Let $\Si$ be a Riemann surface
of positive genus. We embed $X:\Si\to \RR^N$ for some $N\ge 3$ and define 
the center of mass for a function $u\in H^{1,2}(\Si)$ by
$${m_c}(u)=\frac{\int_{\Si}Xe^u}{\int_{\Si}e^u}.$$
Since $\Si$ is of positive genus, we can choose a Jordan curve $\Gamma^1$
on $\Si$ and a closed curve $\Gamma^2$ in $\RR^N\setminus \Si$ such
that $\Gamma^1$ links $\Gamma^2$. We know that
$\inf_{u\in H^{1,2}(\Si)}J_c(u)$ is finite if and only if
$c\in [0,8\pi]$ (see [DJLW]). Now define a family of functions
$h:D\to H^{1,2}(\Si)$ (as in section 2) 
satisfying
$$\lim_{r\to 1}J_{\rho}(h(r,\theta))\to -\infty$$
and
$$ \lim_{r\to 1}{m_c}(h(r,\theta)) \hbox{ as a map from } S^1\to \Gamma^1
\hbox{ is of degree $1$.}$$
Let ${\Cal D}_c$ denote the set of all such families. 
%We  can check that as in Lemma 2.1.
It is also easy to check that ${\Cal D}_c\not=\emptyset$.
Set
$$\alpha_c:=\inf_{h\in {\Cal D}_c}\sup_{u\in h(D)}J_c(u).$$
We first have
$$\alpha_c> -\infty,$$
using the fact that $\Gamma^1$ links $\Gamma^2$ and Lemma 2.2.
Then by the same method as presented above, we can prove that $\alpha_c$
is achieved by some $u_c\in H^{1,2}(\Si)$, which is a solution of
(1.4), for $c\in (8\pi,16\pi)$.
\qed
\enddemo

\enddocument